\documentclass{article}
\usepackage{ijcai22}

\usepackage{times}
\usepackage{soul}
\usepackage{url}
\usepackage[hidelinks]{hyperref}
\usepackage[utf8]{inputenc}
\usepackage[small]{caption}
\usepackage{graphicx}
\usepackage{amsmath}
\usepackage{amsfonts,amssymb}
\usepackage{amsthm}
\usepackage{booktabs}
\usepackage{algorithm}
\usepackage{algorithmic}
\usepackage{makecell}
\usepackage[utf8]{inputenc} % allow utf-8 input
\usepackage[T1]{fontenc}    % use 8-bit T1 fonts
\usepackage{hyperref}       % hyperlinks
\usepackage{url}            % simple URL typesetting
\usepackage{booktabs}       % professional-quality tables
\usepackage{amsfonts}       % blackboard math symbols
\usepackage{nicefrac}       % compact symbols for 1/2, etc.
\usepackage{microtype}      % microtypography
\usepackage{xcolor}         % colors
\usepackage{amsmath}
\usepackage{cleveref}
\usepackage{graphicx}
\usepackage{tabularx}
\usepackage{float}
\usepackage{multirow}
\usepackage{algorithm}
\usepackage{algorithmic}
\usepackage{caption}
\usepackage{subfigure}

\newtheorem{definition}{Definition}

\newcommand{\eat}[1]{}

\title{Learning to Reformulate for Linear Programming}

\author{
Xijun Li$^{1,2}$
Qingyu Qu$^{2}$
Fangzhou Zhu$^{2}$
Jia Zeng$^2$
Mingxuan Yuan$^2$
Kun Mao$^3$\And
Jie Wang$^1$
\affiliations
$^1$MIRA Lab, USTC\\
$^2$Huawei Noah's Ark Lab\\
$^3$Huawei Cloud Co.
\emails

\{xijun.li, quqingyu, zhufangzhou, Zeng.Jia, Yuan.Mingxuan, maokun\}@huawei.com,
jiewangx@ustc.edu.cn
}
\begin{document}

\maketitle

\begin{abstract}
    It has been verified that the linear programming (LP) is able to formulate many real-life optimization problems, which can obtain the optimum by resorting to corresponding solvers such as OptVerse, Gurobi and CPLEX. In the past decades, a serial of traditional operation research algorithms have been proposed to obtain the optimum of a given LP in a fewer solving time. Recently, there is a trend of using machine learning (ML) techniques to improve the performance of above solvers. However, almost no previous work takes advantage of ML techniques to improve the performance of solver from the front end, i.e., the modeling (or formulation). In this paper, we are the first to propose a reinforcement learning-based reformulation method for LP to improve the performance of solving process. Using an open-source solver COIN-OR LP (CLP) as an environment, we implement the proposed method over two public research LP datasets and one large-scale LP dataset collected from practical production planning scenario. The evaluation results suggest that the proposed method can effectively reduce both the solving iteration number (25\%$\downarrow$) and the solving time (15\%$\downarrow$) over above datasets in average, compared to directly solving the original LP instances.
\end{abstract}

\section{Introduction}

% The significance of mathematical programming
Through many years of practices, it has been verified that the mathematical programming (LP)~\cite{2021From} is capable of formulating real-life optimization problems such as planning, scheduling, resource allocation, etc. The LP can obtain the optimum by resorting to corresponding solvers, such as OptVerse~\cite{optverse}, Gurobi~\cite{gurobi}, CPLEX~\cite{cplex}, SCIP~\cite{scip}, etc, which provides the optimal solutions to those real-life optimization problems. Government and business corporation benefit a lot from the practice of mathematical programmings in their daily operations~\cite{mavrotas2021combining}, which thus constantly draws interests from both academics and industry. There are many kinds of mathematical programmings, including linear programming (LP), mixed integer programming (MIP), quadratic programming (QP). In past decades, a collection of classical algorithms (such as simplex~\cite{1987Origins},  barrier~\cite{2003On}, branch and bound~\cite{wolsey2007mixed}, etc.) have been proposed to solve above mathematical programmings and meanwhile they have been implemented and been integrated in above well-established solvers. Amongst these mathematical programmings, LP is the foundation. Thus many performance improvement of solver can be gained from the research of LP.

% The trend of using ml to improve the solver

Recently, there is a trend of using machine learning (ML) techniques to improve traditional combinatorial optimization solvers~\cite{ml4co} on specific problem distributions. Because in real-life scenarios, a practitioner repeatedly solves problem instances from a specific distribution, with redundant patterns and characteristics. For example, managing a large-scale energy distribution network requires solving very similar optimization problems on a daily basis, with a fixed power grid structure while only the demand changes over time. This change of demand is hard to capture by hand-engineered expert rules, and ML-enhanced approaches offer a possible solution to detect typical patterns in the demand history. A serial of machine learning-based approaches have been proposed to improve the performance of above solvers~\cite{khalil2016learning,gasse2019exact,gupta2020hybrid,nair2020solving}.

% Most of work focus on using ml to replace key component of solver and nobody care the formulation
It should be noted that most of above ML-enhanced approaches focus on using ML techniques to replace some key components within solver. Almost no previous work has thought of accelerating the solving of solver from the most front end, i.e., the modeling from a real-life optimization problem to a mathematical programming. Because we unconsciously think that the human experts are totally responsible for modeling and formulation for real-life optimization problems. The expert-designed formulation is deemed as the `perfect' mathematical programming model and sent to the solver to get the optimal solution, which never thinks of how formulation could affect the performance of corresponding solver. But from some optimization theories and empirical studies~\cite{2021From}, the formulation (such as the ordering of variables and constraints) is highly related to both the accuracy and solving speed of solver, which leaves the huge space for improving the performance of solvers through reformulating the mathematical programming.

In this paper, from the perspective of reformulation, we propose a machine learning-based automatic reformulation method for LP, in order to accelerate the solving, where a graph convolutional neural network (GCNN)~\cite{gasse2019exact,nair2020solving} is firstly utilized to capture the patterns and characteristics of variables in the original LP. Then the pattern of variables is sent to a pointer network (PN)~\cite{bello2016neural} from which we can get a new ordering. The new ordering will change the formulation of original LP but still remain its mathematical properties. The parameter of above two neural networks is trained via reinforcement learning (RL). The contributions of this work are summarized as follows:

\begin{itemize}
    \item To our best knowledge, we are the first to propose a machine learning-based reformulation method for linear programming and implement the method using an open source software COIN-OR LP (CLP) as back-end solver.
    \item Extensive experiments have been performed over two public research LP datasets and one large-scale LP dataset collected from practical production planning scenario. The results suggest that the proposed method can effectively reduce both the solving iteration number (25\%$\downarrow$) and the solving time (15\%$\downarrow$) in average, compared to directly solving the original LP instances.
    \item Our proposed method does not restrict the type of solver, which can be implemented over any solver as long as it has corresponding interfaces. Thus the proposed method is a general way to improve the performance of solver.
    
\end{itemize}
\vspace{-3ex}

\section{Background}

% Present the basic knowledge of linear programming

\subsection{Linear programming and its initial basis}

A linear program is an optimization problem of the form

% \[ \mathop{\arg\min}\limits_{\textbf{x}} \left\{ \textbf{c}^T\textbf{x} \mid \textbf{Ax} \leq \textbf{b}, \textbf{l} \leq \textbf{x} \leq \textbf{u} \right\} \]

\begin{equation}
    \mathop{\min}\limits_{\textbf{x}} \textbf{c}^T\textbf{x},
    \quad \textit{s.t.} \quad
    \textbf{Ax} \leq \textbf{b}, \textbf{x} \geq 0
    \label{eq: def_lp}
\end{equation}
where $\textbf{c} \in \mathbb{R}^\textit{n}$ is the objective coefficient vector;
$\textbf{A} \in \mathbb{R}^{\textit{m} \times \textit{n}}$ is the constraint coefficient matrix; and $\textbf{b} \in \mathbb{R}^\textit{m}$ is the constraint right-hand-side vector. $\textbf{x} \in \mathbb{R}^{\textit{n}}$ is the variable vector; 
Usually $\textbf{A}$ is with full row rank. If the condition cannot be met, appropriate unit columns can be added. Note that formulation (\ref{eq: def_lp}) is the standard form of linear programming. Other form of linear programming can convert to the standard form.
% $\textbf{l} $ and $\textbf{u}$ are the lower and upper bounds of variables respectively. 

Suppose that $S$ is a subset of columns of $\textbf{A}$. We use $\textbf{A}_S$ denote the $m \times |S|$ submatrix that contains the columns of $S$. If $S$ is with an order, then the columns of $\textbf{A}_S$ are taken to appear in that order. Similarly, if $d$ is a vector and $S$ is a subset of row indices, then $d_S$ is the corresponding subvector. Again, if $S$ is ordered, the rows of $d_S$ are also subject to that order. With above concepts, the basis can be defined as follows:

\begin{definition}
A basis is a pair $(B, N)$ which splits the variables (columns):
$B=(B_1,...,B_m)\subseteq \{1,...,n\}$ is an ordered subset of column indices such that $\textbf{B}=\textbf{A}_B$ is nonsingular. $B$ is called the basis indices and $\textbf{B}$ is the basis matrix. The variables (columns) $x_j(j\in B)$ are called the basic variables. The remaining variables are nonbasic variables. We use $N=\{j|j\notin B\}$ denote the set of indices of nonbasic variables. 

\end{definition}

Corresponding to each basis, there exists a basic solution $\textbf{x}$ given by
\begin{equation}
    \textbf{x}_N = 0
\end{equation}
\begin{equation}
    \textbf{x}_B = \textbf{B}^{-1}(b-\textbf{A}_N\textbf{x}_N)=\textbf{B}^{-1}b
\end{equation}
The basis $\textbf{B}$ is called feasible if $\textbf{x}_B \geq 0$. The simplex method requires a feasible basis as inputs. If no such basis exists, it usually resorts to solving a auxiliary problem to get a feasible initial basis. Then it continues to solve the original problem with the initial basis. There are three classical methods to construct the initial basis, i.e., all "artificial" basis, the feasible slack and slack basis. Reader of interests can refer to classical textbook~\cite{maros2002computational} on linear programming. Besides, inside commercial solver CPLEX, \cite{bixby1992implementing} proposed a basis construction method that effectively reduces the iteration number of solving over a class of problems. Specifically, they first construct a preference order for all variables including slack variables. The preference order is obtained via comparing the corresponding objective and constraint coefficients of variables. According to the preference order, $m$ variables are selected to construct the basis. However, there is no certain claim what is the best preference order in selection of initial basic variables.

\subsection{Bipartite graph representation of linear programming}

\begin{figure}[t]
	\centering
	\includegraphics[width=0.85\linewidth]{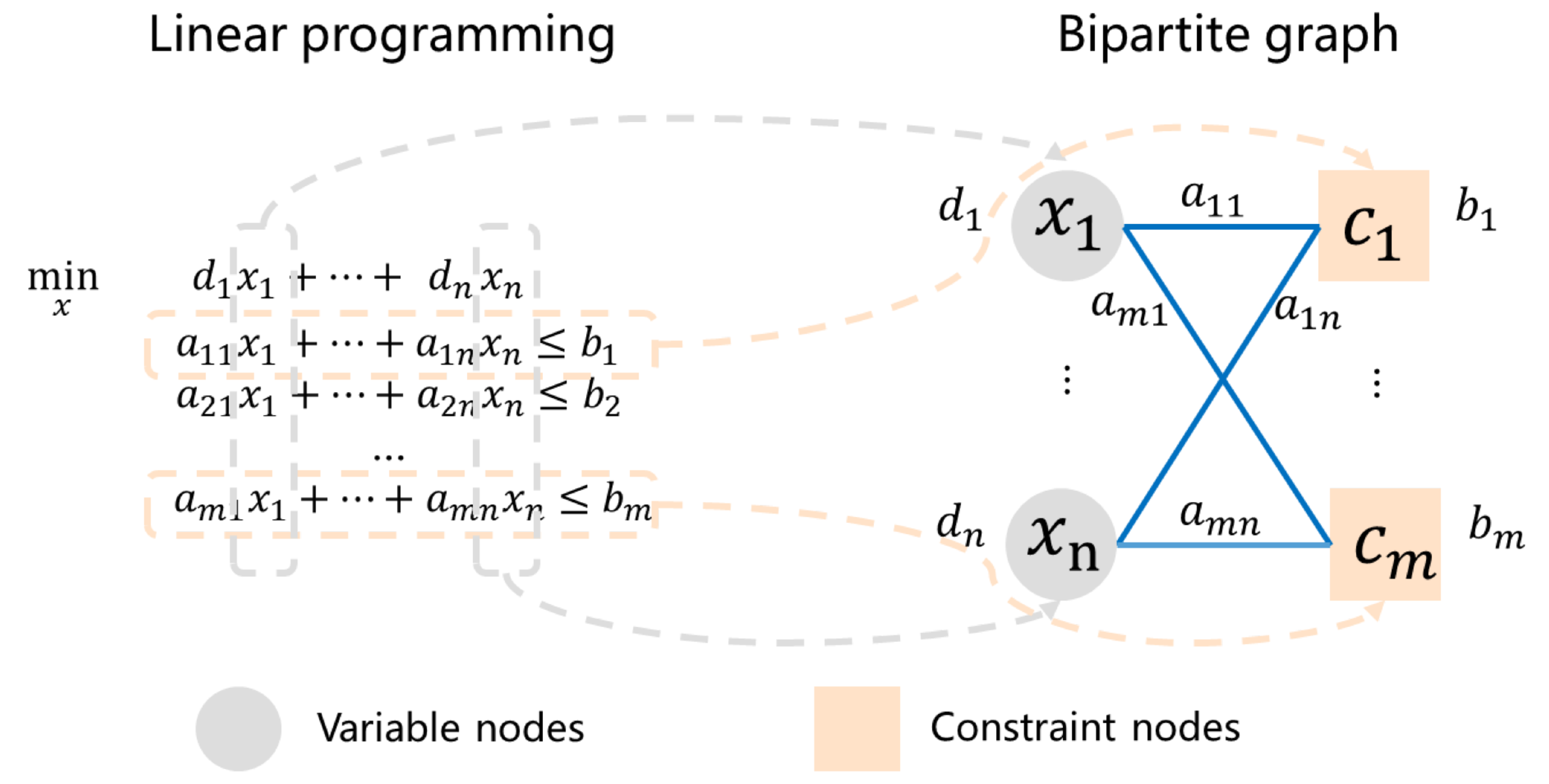}
	\caption{Bipartite graph representation of linear programming}
	\label{fig: bipartite_graph}
\end{figure}

The relation between variables and constraints of LP can be represented by a bipartite graph, where a set of $n$ nodes in the graph represents the $n$ variables contained in the LP and the other set of $m$ nodes correspond to the $m$ constraints of the LP. The edge between a variable node and a constraint node represents the corresponding variable shows in the constraint. The number of edges indicates the number of non-zeros (NNZ) in the constraint matrix $\textbf{A}$. 
An example of the bipartite graph is given in Figure~\ref{fig: bipartite_graph}.
Other information such as the objective coefficients and constraint bounds, etc. can also be added into the bipartite graph. In this way, the lossless representation of the LP can be sent as an input to graph neural networks. Many previous works adopt the representation method or related one to extract high-order embedding information of the original problems, such as Gasses et. al. and Nair et. al. done for mixed integer problem~\cite{nair2020solving,gasse2019exact}, and Selsam et. al. done for Boolean Satisfiability problem~\cite{selsam2018learning}.

\subsection{Pointer network for combinatorial optimization problem}

Combinatorial optimization problems such as Traveling Salesman Problem (TSP), Convex Hull problem, Set Cover problem, etc. play a fundamental role in the development of computer science, which have many applications in manufacturing, planning, genetic engineering, etc. Many kinds of algorithms have been proposed to solve above combinatorial optimization problem, including dynamic programming~\cite{sumita2017improved,chauhan2012survey}, cutting plane~\cite{applegate2003implementing}, local search~\cite{zhang2005novel} and neural network-based search method~\cite{vinyals2015pointer,bello2016neural}. In recent years, the application of neural networks on combinatorial optimization problem has drawn much more attention than other methods~\cite{vinyals2015pointer,bello2016neural}, especially after the proposal of Pointer Network (PN). The pointer network is a sequence-to-sequence model~\cite{vinyals2015pointer} originating from the domain of natural language processing. It can learn the conditional probability of an output sequence of elements that are discrete symbols corresponding to positions in an input sequence, which dedicates to dealing with variable size of output dictionary. Specifically, the pointer network is comprised of two recurrent neural networks (RNN), encoder and decoder. The entire sequence-to-sequence output process is divides as two phases, encoding and decoding. In the encoding phase, the encoder reads initial representation (raw data or linear transformation of raw data) of $s_i $ of input sequence $S=\{s_1, s_2, ..., s_n\}$, one at a time step, and transform them into a sequence of latent memory states $\{enc_i\}_{i=1}^{n}$. In the decoding phase, the decoder network also maintains a latent memory states $\{dec_i\}_{i=1}^{n}$. Then it utilizes the attention mechanism~\cite{vinyals2015pointer} over $\{enc_i\}_{i=1}^{n}$ and $\{dec_i\}_{i=1}^{n}$ to produce a probability distribution over $s_i (i=1,...,n)$. Then one $s_i$ is pointed and output according to the probability distribution and its corresponding decoding embedding $dec_i$ is sent as input to the next decoder step. 
% The entire process is depicted in Figure~\ref{fig: pn}. 
On the TSP, \cite{vinyals2015pointer} trained above neural network in a supervised manner to predict the sequence of visited cities. \cite{bello2016neural} trained the network with reinforcement learning method, using the negative tour length as the reward signal.
\vspace{-2ex}
% \begin{figure}[t]
% 	\centering
    
% 	\includegraphics[width=0.67\linewidth]{figures/pn.pdf}
% 	\caption{The pointer network (PN) architecture}
% 	\label{fig: pn}
% \end{figure}
% \input{content/related_work}
\section{Proposed Solution}

\subsection{Overview}

\begin{figure*}[t]
	\centering
    
	\includegraphics[width=0.85\linewidth]{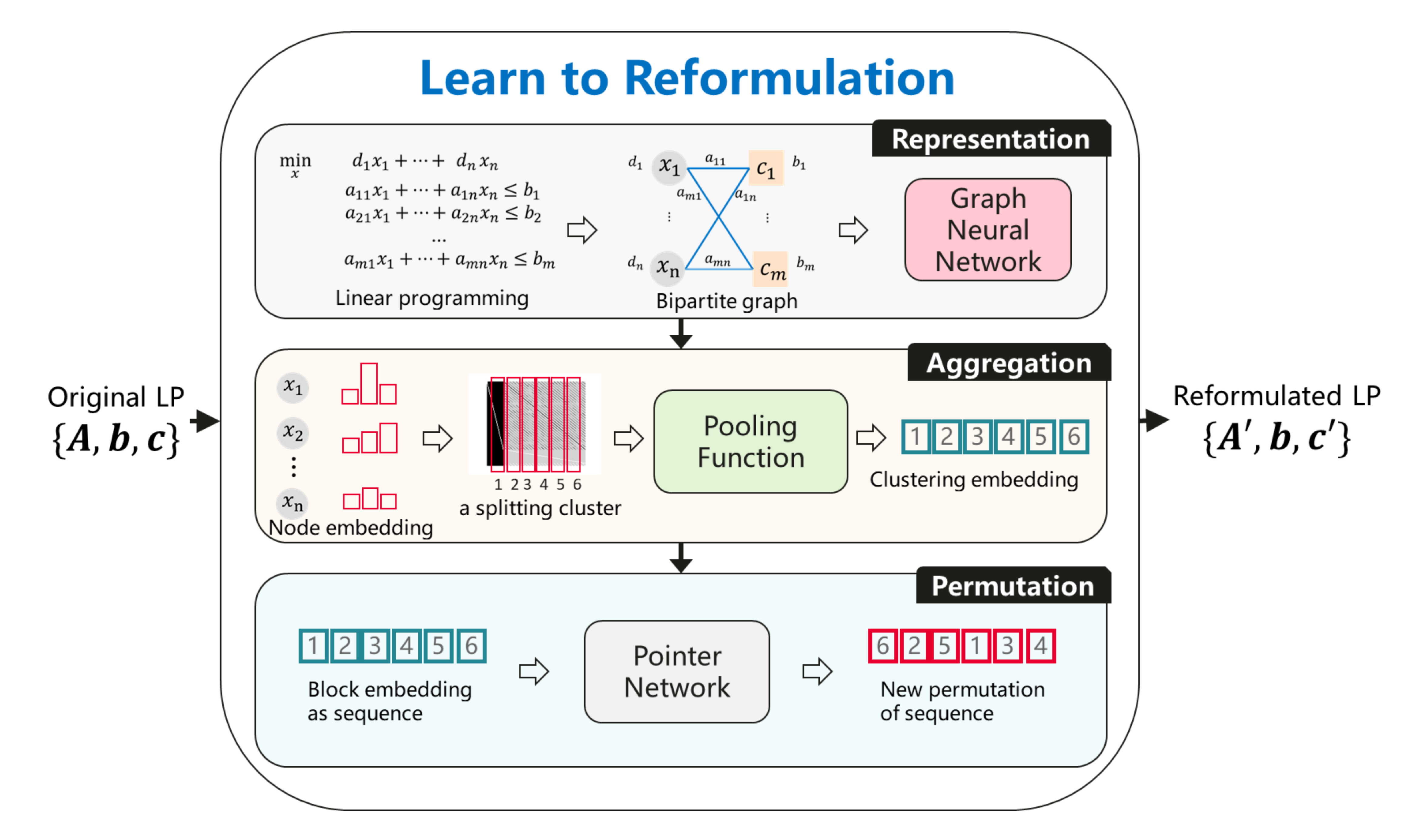}
	\caption{Overview of our proposed machine learning-based automatic reformulation method, which can be summarized in three steps: a) the input LP instance is represented by a bipartite graph and then the embedding of variables of the LP instance is obtained via a graph neural network; b) the embedding of variables will be aggregated with a given splitting cluster of variables and c) taking as input the previous embeddings, a pointer network is used to output a new ordering of variables, i.e., the reformulation of the original LP instance.}
	\label{fig: architecure}
\end{figure*}

In this section, our reformulation method is presented. We firstly introduce how we represent a given linear programming and send it as input into a graph neural network. Then the embedding output by the graph neural network is aggregated with a given variable splittings and passed into a pointer network to get a permutation of variables. The permutation is utilized to reformulate the original linear programming, in order to accelerate the solving process of corresponding solver. The entire process is summarized in Figure~\ref{fig: architecure}.

\subsection{Representation}

We adopt the same method as done in ~\cite{gasse2019exact} to represent a given linear programming as a bipartite graph $(\mathcal{G}, \textbf{C}, \textbf{E}, \textbf{V})$. Specifically, in the bipartite graph, $\textbf{C}\in\mathbb{R}^{m\times c}$ corresponds to the features of the constraints in the LP; $\textbf{V}\in\mathbb{R}^{n\times d}$ denotes the features of the variables in the LP; and an edge $e_{ij}\in \textbf{E}$ between a constraint node $i$ and a variable node $j$ if the corresponding coefficient $A_{i,j} \neq 0$. For simplicity, we just attach the value of $A_{i,j}$ to the corresponding edge $e_{ij}$. Readers of interest can refer to the used features in Appendix.

Next, the bipartite graph representation of LP is sent as input into a two-interleaved graph convolutional neural network (GCNN)~\cite{gasse2019exact}. In detail, the graph convolution is broken into two successive passes, one from variable side to constraint side, and one from constraint side to variable side, which can be formulated as follows:
\begin{equation}
    \textbf{c}_i^{(l+1)} \gets \textbf{f}_{\textbf{C}}\left( \textbf{c}_i^{(l)} , \sum_{j}^{(i,j)\in \textbf{E}} \textbf{g}_{\textbf{C}}(\textbf{c}_i^{(l)}, \textbf{v}_j^{(l)}, e_{ij}) \right), 
\end{equation}
\begin{equation}
    \textbf{v}_j^{(l+1)}  \gets \textbf{f}_{\textbf{V}}\left(\textbf{v}_j^{(l)} , \sum_{i}^{(i,j)\in \textbf{E}} \textbf{g}_{\textbf{V}}(\textbf{c}_i^{(l)}, \textbf{v}_j^{(l)}, e_{ij})\right)
\end{equation}
where $\textbf{f}_{\textbf{C}}, \textbf{g}_{\textbf{C}}, \textbf{f}_{\textbf{V}}$ and $\textbf{g}_{\textbf{V}}$ are 2-layer perceptrons with prenorm layer. We adopt the ReLU as the activation function. And $l$ represents the number of times that we perform the convolution. In our implementation, we set $l=2$. The parameters involved in above GCNN are denoted by $\theta_G$.

\subsection{Aggregation}

The embedding information of variables can be obtained using the GCNN. However, we perform an aggregation operation over the embedding rather than directly sending the them to pointer networks. There are several reasons why we need to perform aggregation. First, the learning capability of pointer network is limited. According to evaluation report of previous work ~\cite{vinyals2015pointer,bello2016neural}, the pointer network can achieve closely optimal results with up to 100 nodes. It performs significantly worse when the number of nodes exceeds 1000. Second, considering all possible permutations of variables of a given LP is intractable. The number of variables of LP that comes from practical scenario can easily exceeds 100. Third, many LPs have its own special structure, which can be exploited to split the variables into several clusters in advance. Many optimization methods such as Benders decomposition~\cite{mavrotas2021combining,gharaei2020joint} have exploited the structure of LP model to accelerate the searching process. Considering all above, we perform the aggregation using the following steps:

\noindent\textbf{Splitting.} For a given LP as shown in~(\ref{eq: def_lp}), the variables $x_i(i=1,...,n)$ are split up into $k$ disjoint clusters $C_j=\{x_{j_1}, x_{j_2}, ..., x_{j_{|C_j|}}\} (j=1,...,k)$. 
Note that the ordering of variables within one cluster is subject to the ordering of variables in the original LP.
The clustering method is not restricted here. It could be specified by human experts or using hyper-graph decomposition method~\cite{manieri2021hyper}.

\noindent\textbf{Pooling function.} With above splitting clusters $C_j (j=1,...,k)$ and variable embedding $\textbf{v}_i (i=1,...,n)$, we perform the aggregation for each cluster via:
\begin{equation}
    \Sigma_{j} = \mathcal{P}(\{\textbf{v}_i|x_{i}\in C_j\})
    \label{eq: pooling}
\end{equation}
where $\mathcal{P}$ is a pooling function which could be maximum, minimum, average or other appropriate functions. We still do not restrict the kind of the pooling function here. $\Sigma_j\in \mathbb{R}^d$ can be understood the embedding of splitting cluster.

\subsection{Permutation}

Given an LP $lp$, we aim to reformulate the LP by reordering the variables of the original LP, in order to improve the solving performance of corresponding solver. 
More specifically, given a sequence of splitting clusters $\{C_j\}_{j=1}^k$ of $lp$, we would like to find a permutation $\pi$ of these clusters to reformulate the original LP. The reformulation is achieved by that 1) between splitting clusters, the variables will be reordered with its cluster according to the permutation $\pi$; and 2) within each cluster, the order of variables remains the same with the original LP. In this way, the coefficients matrix $\mathbf{A}$ and cost coefficients vector $\mathbf{c}$ will correspondingly be altered. We hope the reformulation of original LP can improve the solving performance of solver over the LP. The solving performance can be the solving time, iteration number, solving accuracy (i.e., primal/dual solution violation), etc., which depends on the preference of performance optimization. We formally define the improvement $R$ of solving performance gained from reformulation as:  
\begin{equation}
    R(\pi|lp) = \mathcal{S}_{\mathcal{M}}(lp) - \mathcal{S}_{\mathcal{M}}(lp|\pi)
    \label{eq: reward}
\end{equation}
where $\mathcal{S}_{\mathcal{M}}(lp)$ denotes that with respect to a solving performance metric $\mathcal{M}$, calling a solver to solve $lp$; and $\mathcal{S}_{\mathcal{M}}(lp|\pi)$ refers to that with respect to the same solving metric $\mathcal{M}$, calling the same solver to solve $lp$ reformulated using the variable permutation $\pi$. Using Eq.(\ref{eq: reward}), we can measure how a permutation $\pi$ can improve the solving performance $\mathcal{M}$ of solver, compared to the original LP without reformulation. 

Our aim is to learn a probability distribution $p(\pi|\{C_j\}_{j=1}^k)$ that given a sequence of splitting clusters $\{C_j\}_{j=1}^k$ of $lp$ and corresponding embedding $\{\Sigma_j\}_{j=1}^k$, can assign higher probabilities to "better" permutations and lower probabilities to "worse" ones. The "better" and "worse" are measured using Eq.(\ref{eq: reward}). Similar to~\cite{vinyals2015pointer,bello2016neural}, the probability distribution $p(\pi|\{C_j\}_{j=1}^k)$ utilizes the chain rule to factorize the probability of a permutation as:
\begin{equation}
p(\pi|\{C_j\}_{j=1}^{k}) = \prod_{j=1}^k p(\pi(j)|\pi(<j),\{C_j\}_{j=1}^{k})
\end{equation}
We parameterize $p(\pi|\{C_j\}_{j=1}^k)$ by a pointer network whose parameter is denoted by $\theta_{P}$.

\subsection{Training method}

In our proposed method, there are two main classes of parameters, $\theta_G$ and $\theta_P$, to be learned. Theoretically, the parameters could be trained using supervised learning (SL) as done in~\cite{vinyals2015pointer} or reinforcement learning (RL) as done in~\cite{bello2016neural}. However, we adopt the reinforcement learning method instead of supervised learning method.
First, getting high-quality labelled data (getting improvement $R$ gained from reformulation in our context) is expensive especially when the size of LP is large. Because we need to call a solver to solve two LP instances (i.e., original LP and its reformulated LP) each time we calculate the improvement $R$. Besides, RL is deemed as an effective way to generate better supervised signals, which could help find a more competitive solution than purely supervised learning. Thus we propose to use model-free policy-based RL to learn $\theta_G$ and $\theta_P$.

The training objective is to maximize the expected improvement $R$ over a given LP $lp$, which is formally defined as:
\begin{equation}
    J(\theta_G, \theta_P|lp) = \mathbb{E}_{\pi \sim p_{\theta_G,\theta_P}(.|lp)}[R(\pi|lp)]
\end{equation}
In our training phase, the LPs usually come from the same practical scenario $\mathcal{S}$ such as production planning, bin packing, etc. Thus the total training objective is defined as:
\begin{equation}
    J(\theta_G, \theta_P) = \mathbb{E}_{lp\sim \mathcal{S}}[J(\theta_G, \theta_P|lp)]
    \label{eq: train_objective}
\end{equation}
We utilize stochastic gradient ascent method to optimize Eq.(\ref{eq: train_objective}). According to the REINFORCE algorithm, the gradient of Eq.(\ref{eq: train_objective}) is given as:
\begin{equation}
\begin{aligned}
        \bigtriangledown_{\theta_G,\theta_P}J(\theta_G, \theta_P|lp)=&\mathbb{E}_{\pi \sim p_{\theta_G,\theta_P}(.|lp)}[(R(\pi|lp)-b(lp))\\
        &\bigtriangledown_{\theta_G,\theta_P}\log p_{\theta_G,\theta_P}(\pi|lp)]
\end{aligned}
\end{equation}
where $b(.)$ is a baseline function independent of $\pi$ and estimates the expected improvement $R$ to reduce the learning variances. To enhance the estimate accuracy of baseline function, we additionally introduce a critic network parameterized by $\theta_c$, which is trained with stochastic gradient descent method over a mean squared error (MSE) between the true improvement $R$ and its prediction $b_{\theta_c}(lp)$. The MSE loss is defined as:
\begin{equation}
    \mathcal{L}(\theta_c)=\mathbb{E}_{lp\sim \mathcal{S}, \pi \sim p_{\theta_G,\theta_P}(.|lp)}[b_{\theta_c}(lp)-R(\pi|lp)]^2
\end{equation}
Note that in our implementation, all mentioned-above gradients are approximated using Monte Carlo sampling. We give a snippet of pseudo codes in Algorithm~\ref{alg} (see Appendix).

\section{Experimental Evaluation}

\begin{table*}[t]
    \centering
    \begin{tabular}{cccc}
        \hline
         Dataset & $m$ & $n$ & $NNZ$ \\
         \hline
         BIP & $195.00 \pm 00.00$ & $1083.00 \pm 00.00$ & $7440.00 \pm 00.00$ \\
         WA & $6.43e04 \pm 54.51$ & $6.1e04 \pm 00.00$	& $3.62e05 \pm 6007.41$ \\
         HPP & $1.25e06 \pm 6.93e04$ & $2.66e06 \pm 1.83e05$	 & $6.64e06	\pm 4.29e05$	 \\
      \hline
    \end{tabular}
    \caption{Statistical description of used dataset}
    \label{tab: stat_dataset}
\end{table*}

% Experiments are conducted to verify the effectiveness of the proposed reinforcement learning based reformulation method. Three datasets of LP instances are mainly considered, of which two are public research dataset and one comes from real-world production planning problem. Besides, to get the metric of interests, an open-source LP solver is used as the back-end solver in our experiments.
% \vspace{-1ex}
\subsection{Setting up}
\subsubsection{Implementation detail}
In our implementation, we use the open-source LP solver COIN-OR LP (CLP)~\cite{clp} as the "environment" in the reinforcement learning setting, with which the proposed method interact. To ease the interaction, we first adopt the well-developed CLP python interface library, CyLP~\cite{cylp} as the interface between the proposed reformulation method and CLP solver. Besides, we develop a new interface that can 1) take as input a given permutation; 2) reformulate a LP instance according to the given permutation; and 3) solve the reformulated LP instance and return the solving metric of interest.

With regard to the GCNN and PN involved in the proposed method, we implemented them using PyTorch~\cite{paszke2019pytorch}. The corresponding hyperparameter is summarized in Table~\ref{tab: hyperparameter} (see Appendix). Note that we use the default parameter of CLP solver when we call the solver to solve a given LP.
All experiments are conducted on a computing server, which is equipped with Intel(R) Xeon(R) Platinum 8180M CPU@2.50GHz, a V100 GPU card with 32GB graphic memory and 1TB main memory. 

% \begin{table}
%   \caption{Sample table title}
%   \label{sample-table}
%   \centering
%   \begin{tabular}{lll}
%     \toprule
%     \multicolumn{2}{c}{Part}                   \\
%     \cmidrule(r){1-2}
%     Name     & Description     & Size ($\mu$m) \\
%     \midrule
%     Dendrite & Input terminal  & $\sim$100     \\
%     Axon     & Output terminal & $\sim$10      \\
%     Soma     & Cell body       & up to $10^6$  \\
%     \bottomrule
%   \end{tabular}
% \end{table}
\vspace{-1.5ex}
\subsubsection{Dataset}
The entire evaluation was performed over three sets of Mixed Integer Linear Programming (MILP) problems. All dataset are scenario specific, i.e., they contain problem instances from only a single scenario. Two of them, Balanced Item Placement (BIP) and Workload Apportionment (WA) , are from \cite{ml4co}. The third set is obtained from Huawei Production Planning (HPP). The detailed description about above dataset is summarized in Appendix. Note that we relax the integer constraint of variables in above MILPs to get the LP instance. Within the dataset of BIP and WA, there are respectively in total $10000$ LP instances, which are different in the value of coefficients, the size of constraints and variables. And for the dataset of HPP, there are in total $1000$ LP instances. The statistical description of above datasets is summarized in Table~\ref{tab: stat_dataset}, where $m,n$ and $NNZ$ represent the number of constraints, variables and non-zero coefficients of a linear programming, respectively. For each dataset, $80\%, 10\%$ and $10\%$ of instances are used as the training, validation and testing set respectively. Besides, for each dataset, we train the neural networks involved in the proposed method separately. In other words, we have three sets of $(\theta_G, \theta_P, \theta_c)$ to train.

\begin{figure*}[t]
    \centering
    \subfigure[BIP training data]{
        \includegraphics[width=1.7in]{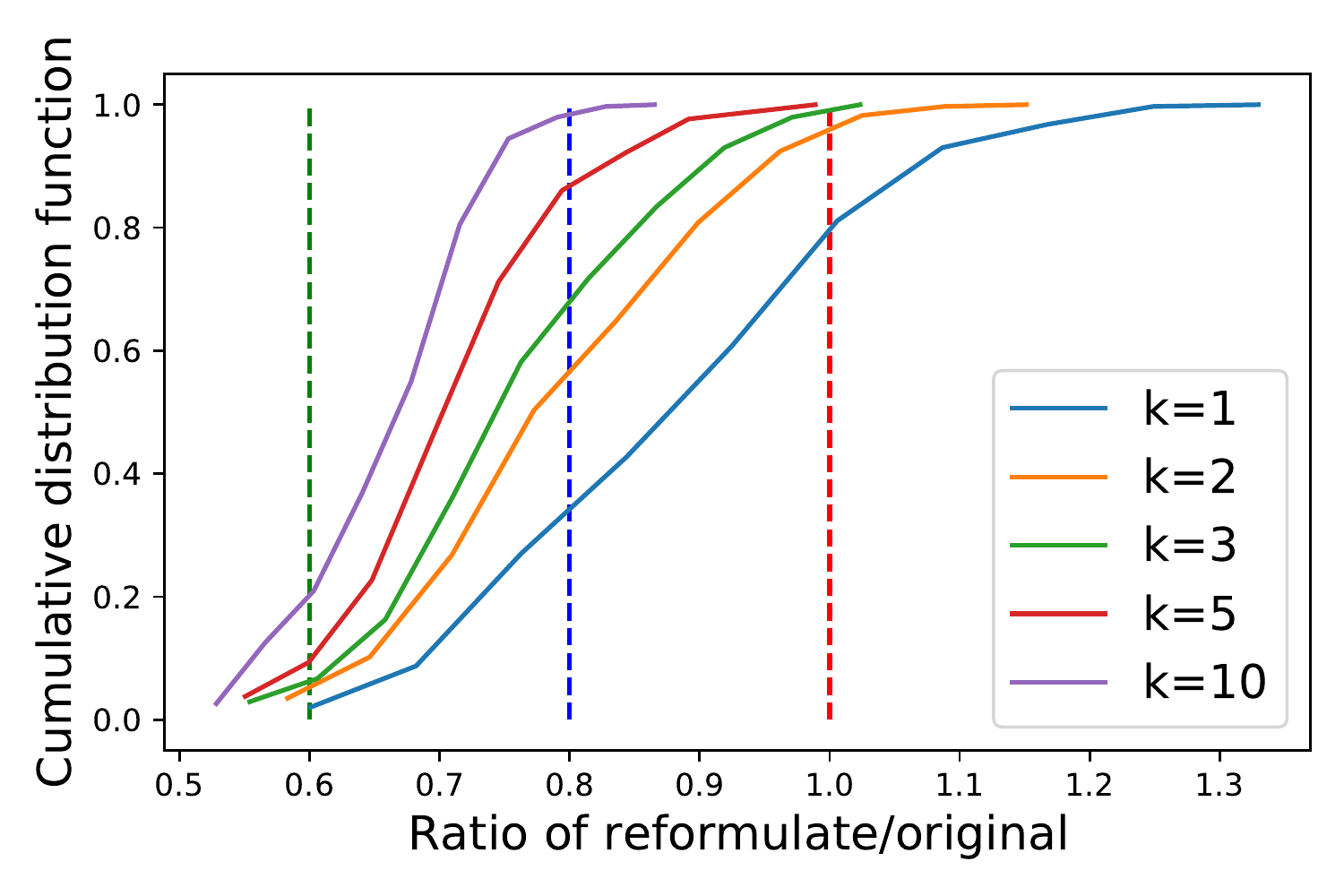}
    }
    \subfigure[WA training data]{
	\includegraphics[width=1.7in]{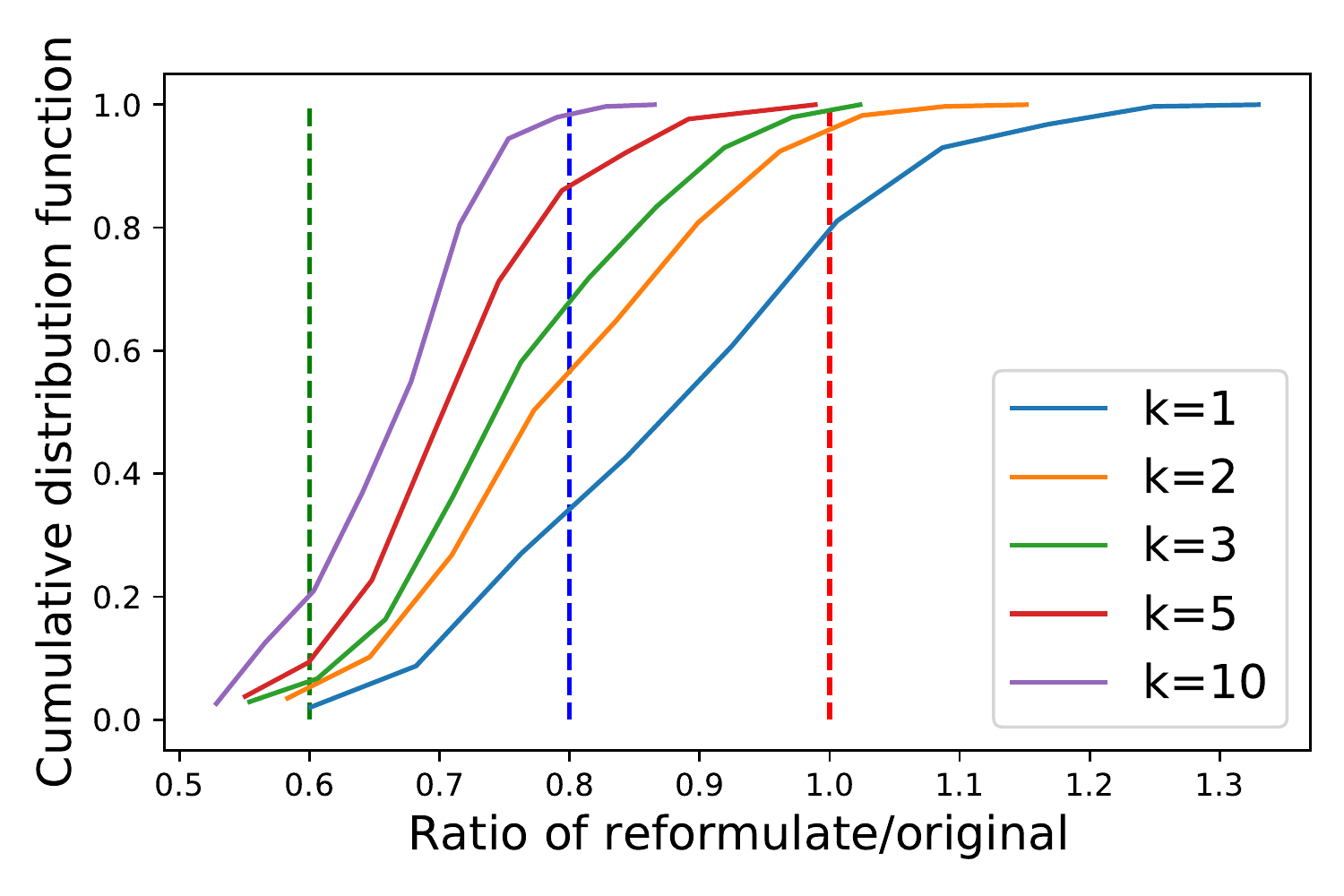}
    }
    \subfigure[HPP training data]{
    	\includegraphics[width=1.7in]{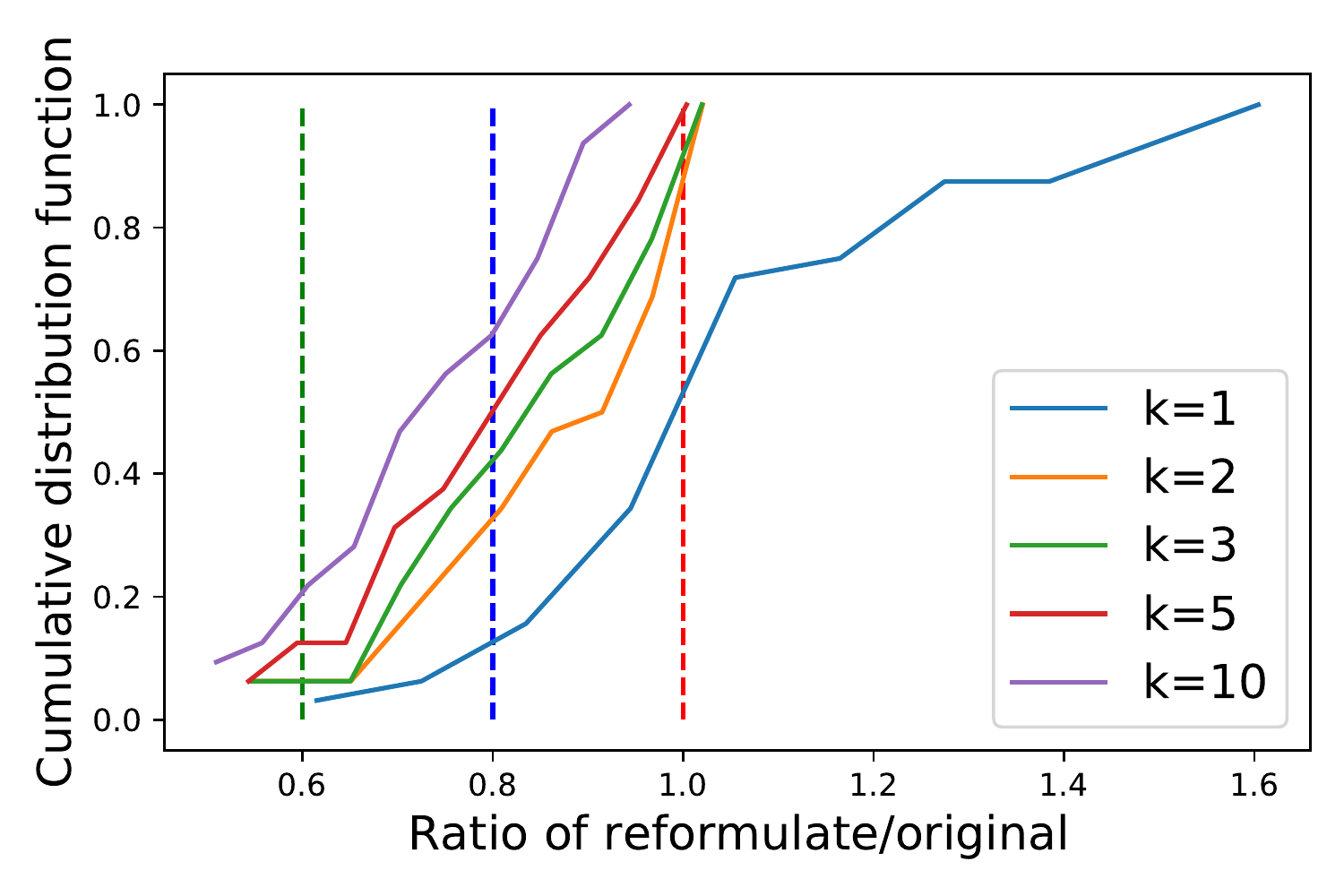}
    }
    \subfigure[BIP testing data]{
        \includegraphics[width=1.7in]{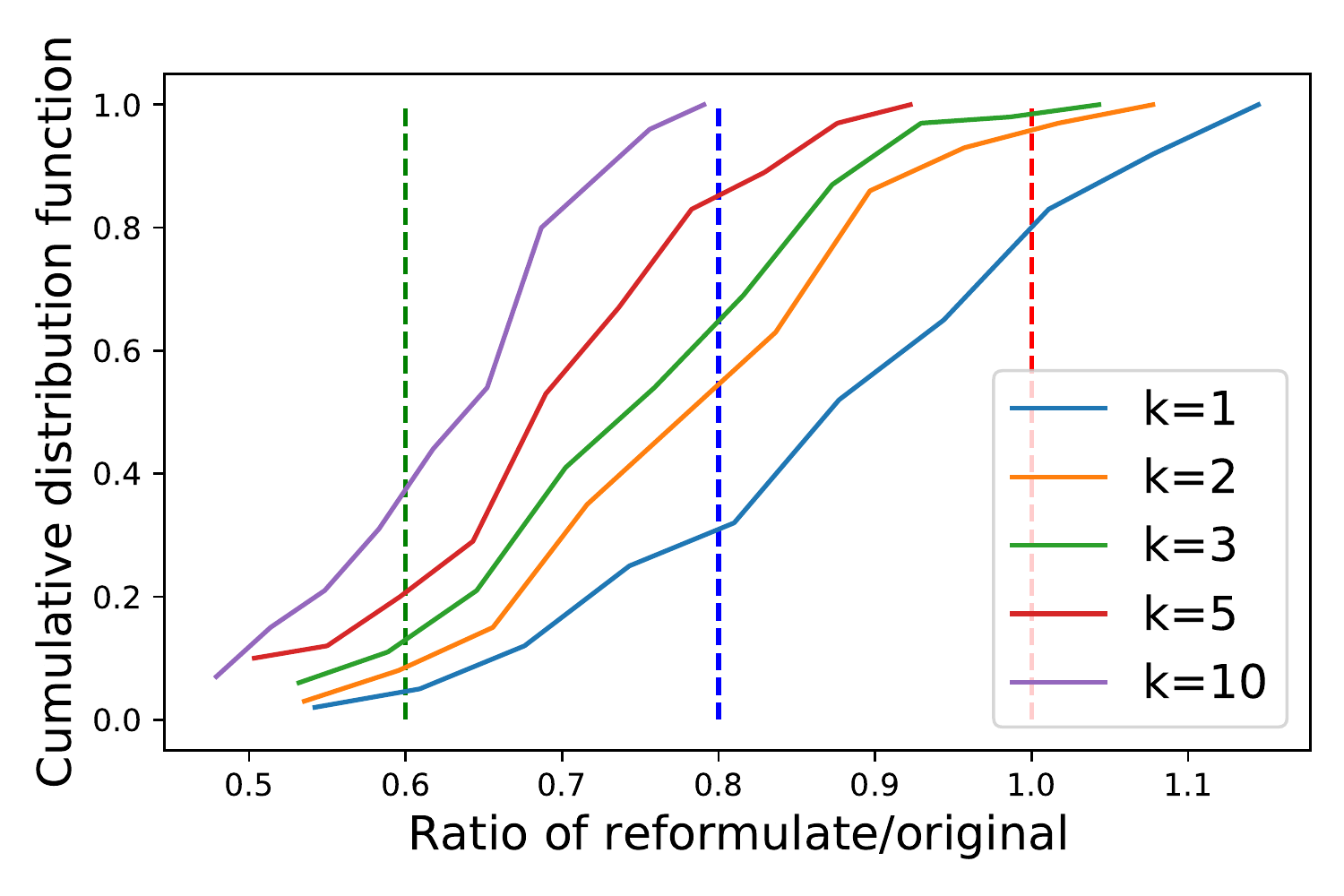}
    }
    \subfigure[WA testing data]{
	\includegraphics[width=1.7in]{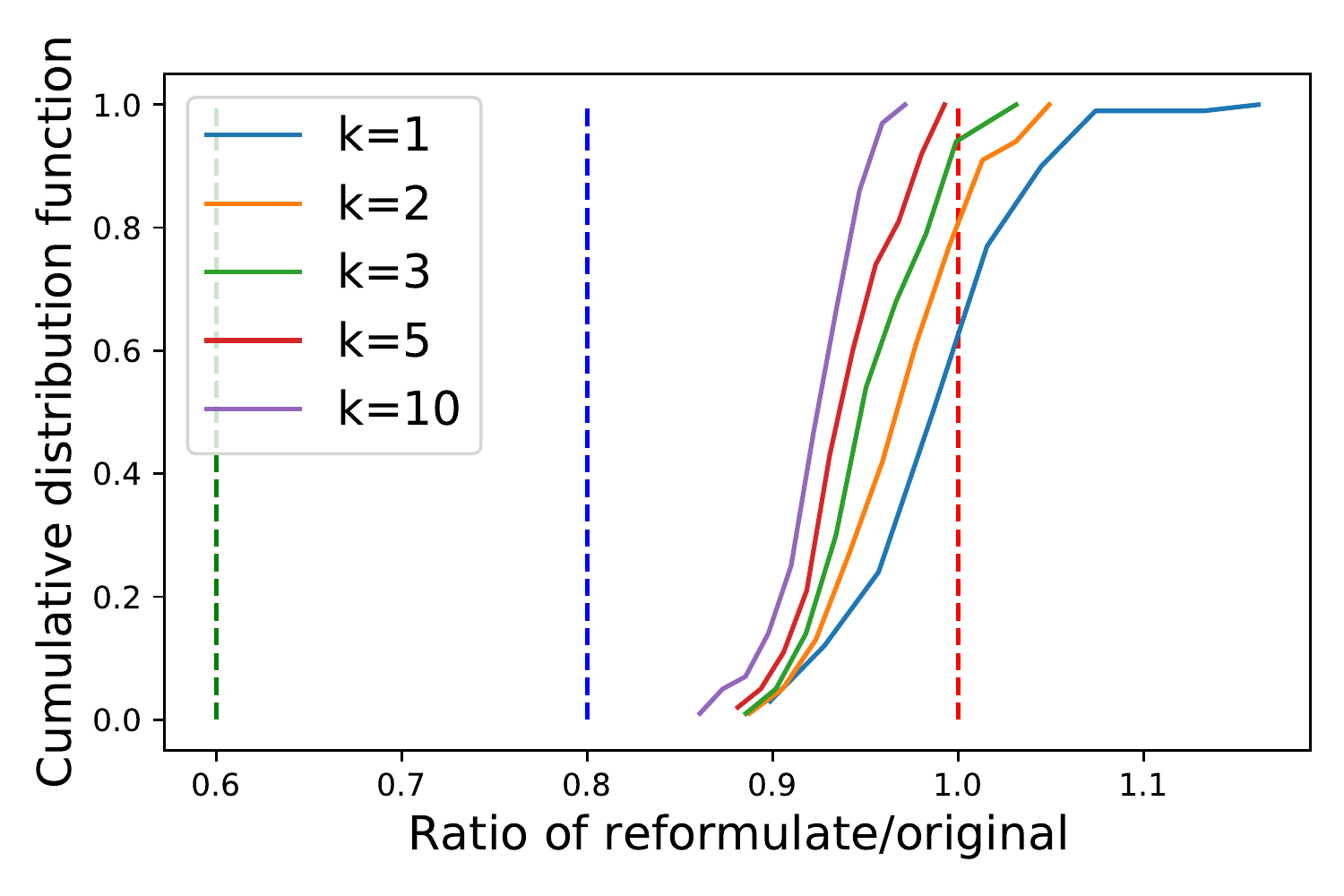}
    }
    \subfigure[HPP testing data]{
    	\includegraphics[width=1.7in]{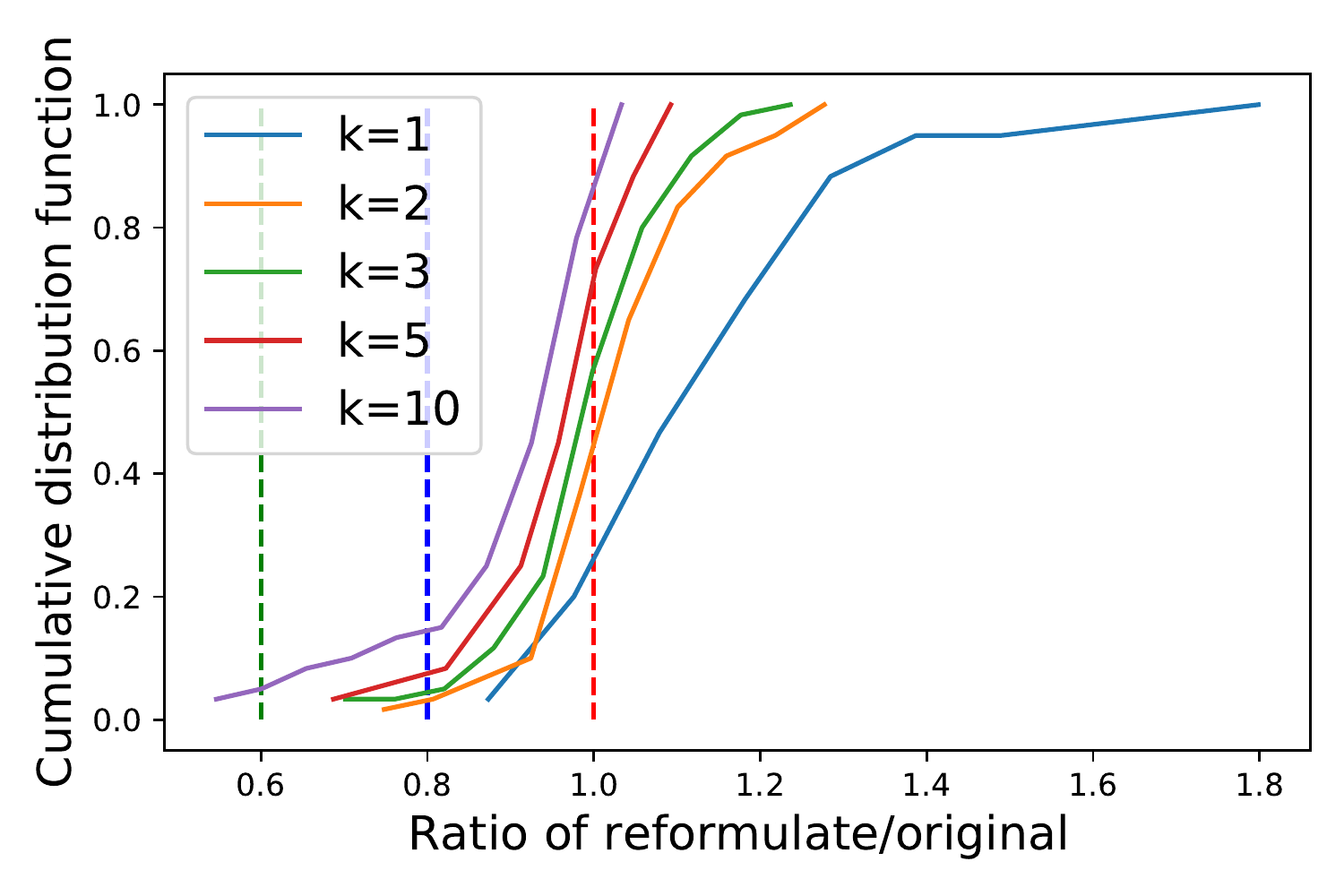}
    }
    \caption{Reduction of the iteration number over training and testing data of three datasets}
    \label{fig: red_iter}
\end{figure*}
\vspace{-1.5ex}
% \begin{figure*}[t]
%     \centering
    
%     \caption{Reduction of iteration numbers over testing dataset}
%     \label{fig: red_iter_test}
% \end{figure*}

% This dataset contains 10000 training instances (pre-split into 9900 train and 100 valid instances).

% This problem deals with apportioning workloads (e.g., data streams) across as few workers (e.g., servers) as possible. The apportionment is required to be robust to any one worker's failure. Each instance problem is modeled as a MILP, using a bin-packing with apportionment formulation. 

% This dataset contains 10000 training instances (pre-split into 9900 train and 100 valid instances).

\subsubsection{Metric of interest}
As defined in Eq.(\ref{eq: reward}), we need to specify the solving performance metric $\mathcal{M}$ when we measure the improvement gained from reformulation. In practice, people usually care about the solving time or iteration number of search process when the solution quality meets a given standard (for example, the maximum of primal/dual infeasible $max\_inf$ of a solution is less than a given primal/dual tolerance). But the solving time might be affected by many factors such as the other background tasks simultaneously running in the same computing server. Thus we finally adopt the iteration number as the solving performance metric $\mathcal{M}$ when we train our model. Besides we also have requirement for solution quality, that is $max\_inf\leq10^{-6}$.

\subsection{Evaluation result}

\subsubsection{Learning convergence}

We first give the learning curves of neural network parameters (i.e., $\theta_P$ and $\theta_G$) in Figure~\ref{fig: convergence} (see Appendix), which shows that the parameter learning of neural networks involved in our method can converge over the three datasets. 
However, compared with small-scale LP instances from BIP, it is relatively harder to learn neural network parameter over the large-scale and complex LP instances from WA and HPP, since the loss variance of former is much smaller than that of latter two.

\subsubsection{Reduction of iteration number}
\label{sec: red_iter_number}
We measure how our proposed reformulation method reduces the iteration number of solving process compared to directly solving the original LP. Specifically, we first call the CLP solver to directly solve LP instances from the testing set of above three datasets and record the iteration number of solving process (here refers to the iteration number of simplex method). Then we reformulate these LP instances using the learned neural networks.
The CLP solver is called again to solve the reformulated LP instances and the iteration number is recorded subsequently. We compare the iteration number of solving reformulated LP instances against that of solving original ones. 
Due to the inference randomness of neural network, here we adopt a $k$-shot inference mechanism, which means that the neural network will infer $k$ times and the best result is kept among the $k$ times of inferences.
The results are presented in Figure~\ref{fig: red_iter}. In these figures, the horizontal axis, Ratio of reformulate/original, denotes the ratio of the iteration number of solving reformulated LP instances to that of solving original ones. The vertical axis denotes the cumulative probability distribution of the ratio over the dataset. 
Several findings can be pointed out: 1) our learning-based reformulation method is effective in reducing the solving iteration number. When $k\geq3$, almost all original LP instance can be solved with fewer solving iteration number via reformulating by our method; 2) our method can be generalized to unseen data since it still performs well over the testing dataset; and 3) On the LP intances from WA and HPP, our method performs slightly worse than on the ones from BIP, which is in accordance with what we observe in the learning convergence.
\vspace{-4ex}
\subsubsection{Reduction of solving time}
We continue to measure how our proposed method reduces the solving time. The same procedure as described in Section~\ref{sec: red_iter_number} is adopted in this experiment. Different from Section~\ref{sec: red_iter_number}, we here compare the solving time between original LP instances and its reformulation obtained from our method, instead of the solving iteration number. Note that for each LP instance reformulated by our method, we keep the best result among $3$ times of inferences of neural networks. The results are presented in Table~\ref{tab: red_solving_time}. From observing the results, several claims can be made: 1) the proposed method can indeed reduce the solving time of given LP instances by reformulating them, which inferably benefits from the reduction of solving iteration number; 2) our method performs slightly worse over the complex and large-scale LP instances but still can reduce at least $8.35\%$ the solving time over the complex LP instances from WA and HPP. Besides, we give a visualization for the reformulation process of our method, in order to figure out what the neural networks have learned (see visual analysis in Appendix A).

\begin{table}[t]
    \centering
    \begin{tabular}{cccc}
        \hline
         Dataset & \quad & Average reduction & Standard error \\ \hline
         \multirow{2}{*}{BIP} 
            & training & $25.18\%$ & $\pm 2.39\%$  \\ \cline{2-4}
            & testing  & $22.35\%$ & $\pm 4.67\%$ \\  \hline
        \multirow{2}{*}{WA} 
            & training & $10.78\%$ & $\pm 0.59\%$  \\ \cline{2-4}
            & testing  & $8.35\%$ & $\pm  1.23\%$ \\  \hline
        \multirow{2}{*}{HPP} 
            & training & $15.62\%$ & $\pm  2.47\%$  \\ \cline{2-4}
            & testing  & $10.47\%$ & $\pm  3.89\%$\\ 
      \hline
    \end{tabular}
    \caption{Reduction of solving time by our proposed method }
    \label{tab: red_solving_time}
\end{table}

\vspace{-2ex}

\section{Conclusion}

In this paper, we propose a machine learning-based reformulation method for LP, in order to improve the solving performance. Specifically, a LP instance is first represented by a bipartite graph, followed by a graph neural network to output the embedding of variables of the LP instance. Then a pointer network takes as input the embedding of variables and output a new ordering of these variables. Then the original LP instance is reformulated according to the new ordering of variables. The parameter of above neural networks is trained using reinforcement learning. Extensive evaluation results over three datasets of LP instances verify the effectiveness of our proposed method in performance improvement of solver.

\clearpage
\bibliographystyle{named}
\bibliography{sample}

\clearpage

\section*{Appendix A}

\subsection*{Feature used in constructing bipartite graph}
\label{sec: features}
% \begin{table}[h]
%     \centering
%     \begin{tabular}{ccc}
%          \hline
%          Tensor & Feature & Description\\
%          \hline
%          \multirow{3}{*}{$\textbf{C}$} & rhs & the right-hand side coefficients of LP, i.e. $\textbf{b}$, normalized with constraint coefficients  \\
%          \cline{2-3}
%          & ub\_cons & upper bound of constraint, normalized with all constraints upper bound \\
%          \cline{2-3}
%          & lb\_cons & lower bound of constraint, normalized with all constraints lower bound \\
%          \hline
%          \multirow{3}{*}{$\textbf{V}$} & obj & the objective coefficients of variables, i.e. $\textbf{c}$  \\
%          \cline{2-3}
%          & lb\_var & upper bound of variable, normalized with all variables lower bound \\
%          \cline{2-3}
%          & ub\_var & lower bound of variable, normalized with all variables upper bound \\
%          \hline
%          $\textbf{E}$ & coef & the constraint coefficients of variables, i.e. $\textbf{A}$ , normalized per constraint \\
%          \hline
%     \end{tabular}
%     \caption{Used features of constraints, variables and edges in the bipartite graph}
%     \label{tab: features}
% \end{table}

\noindent{\textbf{Features of constraint:}}
\begin{itemize}
    \item rhs: the right-hand side coefficients of LP, i.e. $\textbf{b}$, normalized with constraint coefficients.
    \item ub\_cons: upper bound of constraint, normalized with all constraints upper bound.
    \item lb\_cons: lower bound of constraint, normalized with all constraints lower bound.
\end{itemize}

\noindent{\textbf{Features of variable:}}
\begin{itemize}
    \item lb\_var: upper bound of variable, normalized with all variables lower bound.
    \item lb\_var: upper bound of variable, normalized with all variables lower bound.
    \item ub\_var: lower bound of variable, normalized with all variables upper bound.
\end{itemize}

\noindent{\textbf{Features of edge:}}
\begin{itemize}
    \item coef: the constraint coefficients of variables, i.e. $\textbf{A}$, normalized per constraint
\end{itemize}

\subsection*{Pseudo code of training method}
\begin{algorithm}
	\caption{Training procedure of the proposed method} 
	\label{alg} 
	\begin{algorithmic}
		\REQUIRE
		\begin{itemize}
		    \item $\mathcal{S}$: set of linear programming problems
		    \item $ T$: number of training steps
		    \item $ B $: batch size
		\end{itemize}
% 		\REQUIRE $\mathcal{S}$: set of linear programming problems
% 		\REQUIRE $ T$ : number of training steps
% 		\REQUIRE $ B $: batch size
		\ENSURE
		\begin{itemize}
		    \item $\theta_G$: parameters for graph convolutional neural network
		    \item $\theta_P $: parameters for pointer network
		    \item $\theta_c $: parameters for critic network
		\end{itemize}

		\FOR{steps $ t = 1 $ to $T$}
		\STATE $ lp_i \sim \mathcal{S}$ for $i\in\{1,...,B\}$ 
		\STATE $\pi_i \sim p_{\theta_G,\theta_P}(.|lp_i)$ for $i\in\{1,...,B\}$
		\STATE Calculate $R_i$ using Eq.(\ref{eq: reward}) for $i\in\{1,...,B\}$
		\STATE $b_i=b_{\theta_c}(lp_i)$ for $i\in\{1,...,B\}$
		\STATE $g \gets \frac{1}{B}\sum_{i=1}^{B}(R_i-b_i)\bigtriangledown_{\theta_G,\theta_P}\log p_{\theta_G,\theta_P}(\pi_i|lp_i)$
		\STATE $ \mathcal{L}_c \gets \frac{1}{B}\sum_{i=1}^{B}(b_i-R_i)^2$
		\STATE Perform a gradient ascent step to update $ \theta_G $ and $\theta_P$ using $g$ respectively
		\STATE Perform a gradient descent step to update $ \theta_c $ using $\bigtriangledown_{\theta_c}\mathcal{L}_c$
		\ENDFOR
		\STATE return $\theta_G, \theta_P$ and $\theta_c$
	\end{algorithmic} 
\end{algorithm}

\subsection*{Hypeparameters setting}

Part of important hyperparameters involved in our method is list in Table~\ref{tab: hyperparameter}.

\begin{table}[h]
    \centering
    \begin{tabular}{cc}
        \hline
         Name & Used value\\
         \hline
         Optimizer& \textit{ADAM}  \\
         \# epochs ($T$) & 40000 \\
         \# splitting cluster & 20 \\
         
         Batch size ($B$) & 8 \\
         Train size & 640 \\
         Validation size & 320 \\
         Learning rate & $10^{-4}$ \\
         Decay ratio of learning rate & 0.96 \\
         Gradient clip normalizer & $l_2$ norms \\
         Dimension of input embedding in PN& 128 \\
         Dimension of hidden layers in PN & 128 \\
         Dimension of input embedding in GCNN & 64 \\
         \# of times that performs convolution ($l$) & 2\\
      \hline
    \end{tabular}
    \caption{Hyperparameters setting}
    \label{tab: hyperparameter}
\end{table}

\subsection*{Description of dataset}

\noindent\textbf{Balanced Item Placement.}
This problem deals with spreading items (e.g., files or processes) across containers (e.g., disks or machines) utilizing them evenly. Items can have multiple copies, but at most, one copy can be placed in a single bin. The number of items that can be moved is constrained, modeling the real-life situation of a live system for which some placement already exists. Each problem instance is modeled as a MILP, using a multi-dimensional multi-knapsack formulation. 

\noindent\textbf{Workload Apportionment.}
This problem deals with apportioning workloads (e.g., data streams) across as few workers (e.g., servers) as possible. The apportionment is required to be robust to any one worker's failure. Each instance problem is modeled as a MILP, using a bin-packing with apportionment formulation. 

\noindent\textbf{Huawei Production Planning.}
The planning and scheduling optimization problems solved in the Huawei production planning engine.

\subsection*{Visual analysis}

\begin{figure*}[t]
    \centering
    \subfigure[original LP1]{
        \includegraphics[width=1.2in]{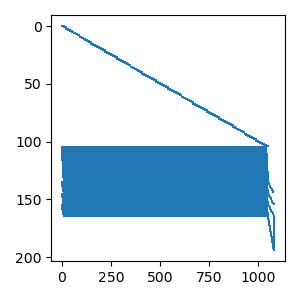}
    }
    \subfigure[reformulated LP1]{
	\includegraphics[width=1.2in]{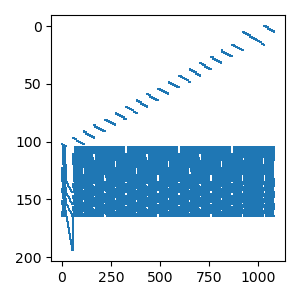}
    }
    \subfigure[original LP2]{
    	\includegraphics[width=1.2in]{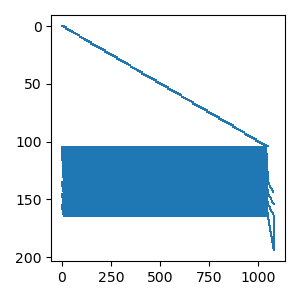}
    }
    \subfigure[reformulated LP2]{
    	\includegraphics[width=1.2in]{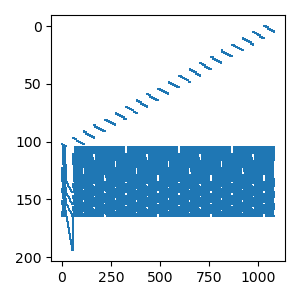}
    }
    \caption{Visualization of reformulation process over BIP dataset}
    \label{fig: vis_matrix1}
\end{figure*}

\begin{figure*}[t]
    \centering
    \subfigure[original LP1]{
        \includegraphics[width=1.2in]{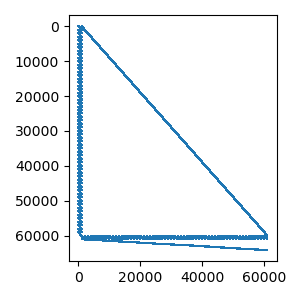}
    }
    \subfigure[reformulated LP1]{
	\includegraphics[width=1.2in]{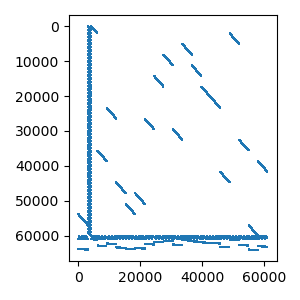}
    }
    \subfigure[original LP2]{
    	\includegraphics[width=1.2in]{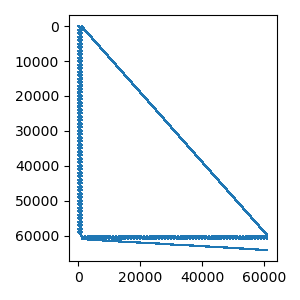}
    }
    \subfigure[reformulated LP2]{
    	\includegraphics[width=1.2in]{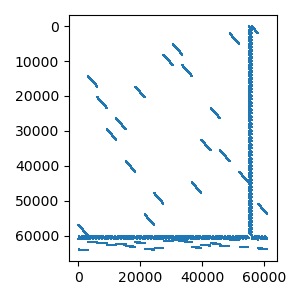}
    }
    \caption{Visualization of reformulation process over WA dataset}
    \label{fig: vis_matrix2}
\end{figure*}

\begin{figure*}[t]
    \centering
    \subfigure[original LP1]{
        \includegraphics[width=1.2in]{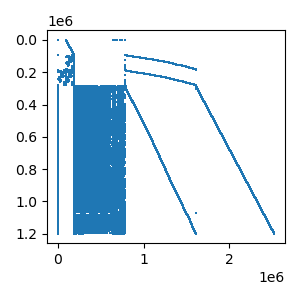}
    }
    \subfigure[reformulated LP1]{
	\includegraphics[width=1.2in]{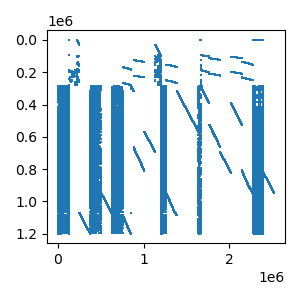}
    }
    \subfigure[original LP2]{
    	\includegraphics[width=1.2in]{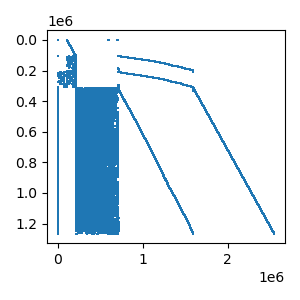}
    }
    \subfigure[reformulated LP2]{
    	\includegraphics[width=1.2in]{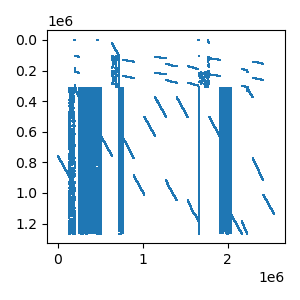}
    }
    \caption{Visualization of reformulation process over HPP dataset}
    \label{fig: vis_matrix3}
\end{figure*}

In order to figure out what the neural network have learned, we give the visualization for the reformulation process of our method, over the three datasets of LP instances. Specifically, from each dataset, we select two LP instances. We reformulate them by the learned neural network of our proposed method. Then we draw the coefficient matrix of original LP instances and reformulated ones respectively, which 
are shown in Figure~\ref{fig: vis_matrix1} to Figure~\ref{fig: vis_matrix3}. Observing these figures, we can find that 1) our proposed reformulation method indeed capture the characteristics of LP instances originated from different scenarios. Because the pattern of corresponding reformulated LP instances are greatly different across different datasets but are similar between LP instances within the same dataset; 2) The reformulation is relatively stable when the original LP instances are highly similar. All the original LP instances of BIP are with the same number of constraints and variables, which is only different in the value of coefficients. Thus the pattern of the corresponding reformulated LP instances are almost the same (see Figure~\ref{fig: vis_matrix1}). However, the pattern of reformulated LP instances of WA and HPP is quiet different (see Figure~\ref{fig: vis_matrix2} and Figure~\ref{fig: vis_matrix3}) because the corresponding original LP instances differs in not only the value of coefficients but also the number of constraints and variables.

\clearpage
\begin{figure*}[t]
    \centering
    \subfigure[loss curve of $\theta_P$ over BIP]{
        \includegraphics[width=1.7in]{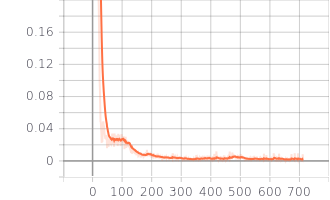}
    }
    \subfigure[loss curve of $\theta_P$ over WA]{
	\includegraphics[width=1.7in]{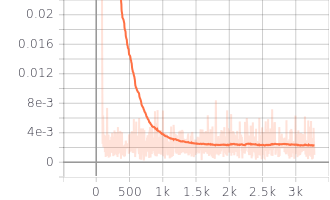}
    }
    \subfigure[loss curve of $\theta_P$ over HPP]{
    	\includegraphics[width=1.7in]{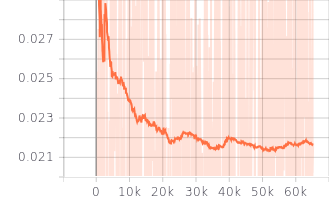}
    }
    
    \subfigure[loss curve of $\theta_G$ over BIP]{
        \includegraphics[width=1.7in]{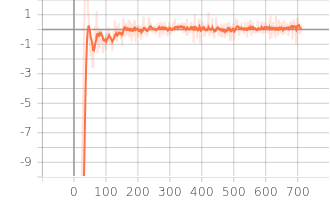}
    }
    \subfigure[loss curve of $\theta_G$ over WA]{
	\includegraphics[width=1.7in]{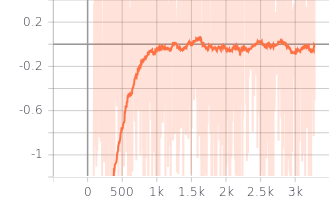}
    }
    \subfigure[loss curve of $\theta_G$ over HPP]{
    	\includegraphics[width=1.7in]{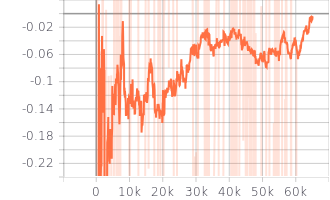}
    }
    \caption{Learning convergence result of our proposed method over above three datasets}
    \label{fig: convergence}
\end{figure*}

% \clearpage
% \input{content/appendix_c}

\end{document}